\documentclass[12pt]{article}
\usepackage{amsmath, graphicx, amsfonts,amssymb}
\usepackage{color}

\def\be{\begin{equation}}
\def\ee{\end{equation}}
\def\bea{\begin{eqnarray}}
\def\eea{\end{eqnarray}}
\def\bt{\begin{theorem}}
\def\et{\end{theorem}}
\def\bl{\begin{lemma}}
\def\el{\end{lemma}}
\def\br{\begin{remark}}
\def\er{\end{remark}}
\def\bc{\begin{corollary}}
\def\ec{\end{corollary}}
\def\bd{\begin{definition}}
\def\ed{\end{definition}}
\def\a{\alpha}
\def\b{\beta}

\def\d{\delta}

\def\k{\kappa}

\def\m{\mu}

\def\s{\sigma}

\def\bbR{\mathbb{R}}

\def\dis{\displaystyle}

\def\bn{B_2^}
\def\cn{B_{\infty}^}
\def\b1{B_{1}^}

\def\pt{\partial}

\def\ba{\begin{array}}
\def\ea{\end{array}}
\def\ben{\begin{enumerate}}
\def\een{\end{enumerate}}
\newtheorem{theorem}{Theorem}[section]
\newtheorem{lemma}{Lemma}[section]
\newtheorem{remark}{Remark}[section]
\newtheorem{proposition}{Proposition}[section]
\newtheorem{corollary}{Corollary}[section]

\newtheorem{definition}{Definition}[section]
\begin{document}
\title{New $L_p$ Affine Isoperimetric Inequalities
}

\author{Elisabeth Werner \thanks{Partially supported by an NSF grant, a FRG-NSF grant and  a BSF grant}
\and Deping Ye }
\date{}
\maketitle
\begin{abstract}
We prove new $L_p$ affine isoperimetric inequalities
for all $ p \in [-\infty,1)$.  We establish, for all $p\neq -n$, a duality
formula  which shows that $L_p$ affine surface area of a convex
body $K$ equals $L_\frac{n^2}{p}$ affine surface area of the polar
body $K^\circ$.
\end{abstract}

\section{Introduction}
An affine isoperimetric inequality relates two functionals
associated with convex bodies (or more general sets) where the
ratio of the functionals is invariant under non-degenerate linear
transformations. These affine isoperimetric inequalities are more
powerful than their better known Euclidean relatives.

\par
This article deals with affine isoperimetric inequalities for the $L_p$ affine surface area.
$L_p$ affine surface area was introduced by Lutwak in the ground
breaking paper \cite{Lu1}. It is now at the core of
the rapidly developing  $L_p$ Brunn Minkowski theory.
Contributions here include new interpretations of $L_p$ affine
surface areas \cite{MW2, SW4, SW5}, the discovery of new
ellipsoids \cite{Lud2, LuYZh1}, the study of solutions of
nontrivial ordinary and, respectively, partial differential
equations (see e.g. Chen \cite{CW1}, Chou and Wang \cite{CW2},
Stancu \cite{SA1, SA2}), the study of the $L_p$ Christoffel-Minkowski problem by
Hu, Ma and Shen \cite{HMS},  a new proof by Fleury, Gu\'edon and
Paouris \cite{FGP}  of a result by Klartag \cite{Kl} on
concentration of volume, and characterization theorems by Ludwig
and Reitzner \cite{LR1}.
\par
The case $p=1$ is the classical affine surface area which goes
back to Blaschke \cite{Bl1}. Originally a basic affine invariant
from the field of affine differential geometry, it has recently
attracted increased attention too (e.g. \cite{Bar,
lei,
Lu2, MW1, SW1}). It  is fundamental  in the theory of valuations
(see e.g.  \cite{A1, A2, LudR, K}), in approximation of convex bodies by
polytopes \cite{Gr2, SW5, LuSchW} and it is the subject of the
affine Plateau problem solved in ${\mathbb R}^3$ by Trudinger and
Wang \cite{TW1, Wa}.

The classical affine isoperimetric inequality which gives an upper
bound for the affine surface area in terms of volume proved to be
the key ingredient in many problems (e.g. \cite{Ga, GaZ, Lu-O,
Sch}). In particular, it was used to show the uniqueness of
self-similar solutions of the affine curvature flow and to study
its asymptotic behavior by Andrews \cite{AN1, AN2}, Sapiro and
Tannenbaum \cite{ST1}.

$L_p$ affine isoperimetric inequalities were first established by
Lutwak for $p >1$ in \cite{Lu1}.
 There has been  a growing body of work in this
area since from which we quote only Lutwak, Yang and Zhang
\cite{LYZ1, LYZ2} and Campi and Gronchi \cite{CG}.
 \par
Here we derive new $L_p$ affine isoperimetric inequalities for all
$ p \in [-\infty,1)$. We give new interpretations of $L_p$ affine
surface areas. We establish, for all $p\neq -n$, a duality formula
which shows that $L_p$ affine surface area of a convex body $K$
equals $L_\frac{n^2}{p}$ affine surface area of the polar body
$K^\circ$. This formula was proved in \cite{Hug} for $p>0$. \vskip
5mm From now on we will always assume that  the centroid of a
convex body $K$ in $\mathbb R^n$ is at the origin. We  write $K\in
C^2_+$ if $K $ has $C^2$ boundary with everywhere strictly
positive Gaussian curvature. For real  $p \neq -n$, we define  the
$L_p$ affine surface area $as_{p}(K)$ of $K$ as in \cite{Lu1} ($p
>1$) and \cite{SW5} ($p <1$) by
\begin{equation} \label{def:paffine}
as_{p}(K)=\int_{\partial K}\frac{\kappa_K(x)^{\frac{p}{n+p}}}
{\langle x,N_{ K}(x)\rangle ^{\frac{n(p-1)}{n+p}}} d\mu_{ K}(x)
\end{equation}
and
\begin{equation}\label{def:infty}
as_{\pm\infty}(K)=\int_{\partial K}\frac{\kappa _K (x)}{\langle
x,N_{K} (x)\rangle ^{n}} d\mu_{K}(x)
\end{equation}
provided the above integrals exist. $N_K(x)$ is the outer unit
normal vector at $x$ to $\partial K$, the boundary of $K$. $\k
_K(x)$ is the Gaussian curvature at $x\in \pt K$ and $\m _K$
denotes the usual surface area measure on $\pt K$.  $\langle
\cdot, \cdot\rangle$ is the standard inner product on $\bbR^n$
which induces the Euclidian norm $\|\cdot\|$.
In particular, for $p=0$
$$
as_{0}(K)=\int_{\partial K} \langle x,N_{ K}(x)\rangle
\,d\mu_{K}(x) = n|K|,
$$
where $|K|$ stands for the $n$-dimensional volume of $K$. More generally, for a set $M$, $|M|$
denotes the Hausdorff content of its appropriate dimension.
For $p=1$
$$
as_{1}(K)=\int_{\partial K}\kappa_K(x)^{\frac{1}{n+1}} d\mu_{ K}(x)
$$
is the classical affine surface area  which is independent
of the position of $K$ in space.

If the boundary of $K$ is sufficiently smooth then (\ref{def:paffine})
and (\ref{def:infty}) can be written as integrals over the boundary
$\pt B^n_2=S^{n-1}$ of the Euclidean unit ball $B^n_2$  in $\mathbb R^n$
$$
as_{p}(K)=\int_{S^{n-1}}\frac{f_{K}(u)^{\frac{n}{n+p}}}
{h_K(u)^{\frac{n(p-1)}{n+p}}}
d\sigma(u).
$$
$\s$ is the usual surface area measure on $S^{n-1}$.
$h_{K}(u)$ is the support function of direction $u\in S^{n-1}$,
and $f_{K}(u)$ is the curvature function, i.e. the reciprocal of
the Gaussian curvature $\kappa _K(x)$ at this point $x \in
\partial K$ that has $u$ as outer normal. In particular, for
$p=\pm \infty$,
\begin{equation}\label{inf-aff}
as_{\pm\infty}(K)
=\int_{S^{n-1}}\frac{1}{h_K(u)^{n}}
d\sigma(u)
=n|K^{\circ}|
\end{equation}
where $K^\circ=\{y\in \bbR^n, \langle x, y\rangle\leq 1, \forall
x\in K\}$ is the polar body of $K$. \vskip 4mm In Sections 2 and 3
we give new geometric interpretations of the $L_p$ affine surface
areas and obtain as a consequence \vskip 3mm \noindent {\bf
Corollary \ref{general-p-affine-surface1}} {\em Let $K$ be a
convex body in $C^2_+$ and let $p\neq -n$ be a real number. Then
$$
as_{p}(K)=as_{\frac{n^2}{p}}(K^\circ).
$$
}
In Section 4  we prove  the  following new $L_p$ affine isoperimetric inequalities.
 For $p \geq 1$ they were proved by Lutwak  \cite{Lu1}.
 \vskip 3mm
\noindent
{\bf Theorem \ref{p-aff-iso4}} {\em Let $K$ be a convex
body with centroid at the origin.
\par
(i)
If $p\geq 0$,  then
\begin{eqnarray*}
\frac{as_p(K)}{as_p(B^n_2)}\leq \left(\frac{|K|}{|\bn
n|}\right)^{\frac{n-p}{n+p}},
\end{eqnarray*}
with equality if and only if $K$ is an ellipsoid. For $p=0$, equality holds trivially for all $K$.
\par
(ii) If $-n<p<0$, then

\begin{eqnarray*}
\frac{as_p(K)}{as_p(B^n_2)}\geq
\left(\frac{|K|}{|B^n_2|}\right)^{\frac{n-p}{n+p}},
\end{eqnarray*} with
equality if and only if $K$ is an ellipsoid.
\par
(iii) If  $K$ is in addition  in $C^2_+$ and if $p < -n$, then
\begin{equation*}
c^{\frac{np}{n+p}}\left(\frac{|K|}{|\bn
n|}\right)^{\frac{n-p}{n+p}} \leq \frac{as_p(K)}{as_p(B^n_2 )}.
\end{equation*}}The constant $c$ in (iii) is the constant from the Inverse
Santal\'o inequality due to Bourgain and Milman \cite{BM}. This
constant has recently been improved by Kuperberg \cite{GK2}.
We give examples  that the above isoperimetric inequalities cannot be
improved.

In Theorem \ref{p-aff-sio8}  we show a monotonicity behavior of
the  quotient $\left(\frac{as_r(K)}{n|K|}\right)^\frac{n+r}{r}$,
namely
$$
\bigg(\frac{as_r(K)}{n|K|}\bigg)\leq
\bigg(\frac{as_t(K)}{n|K|}\bigg)^{\frac{r(n+t)}{{t(n+r)}}}.
$$
and as a consequence obtain

\vskip 3mm
\noindent
{\bf Corollary \ref{p-aff-iso3}}
{\em
Let
$K$ be convex body in $\mathbb R^n$ with centroid at the origin.

(i) For all $p\geq 0$
\begin{equation*}
as_p(K)\ as_p(K^\circ)\leq n^2|K|\ |K^\circ|.
\end{equation*}

(ii) For $-n<p<0$, $$ as_p(K)as_p(K^\circ)\geq n^2|K|\ |K^\circ|.
$$
If $K$ is in addition in $C^2_+$,  inequality  (ii) holds for all
$p<-n$. }

%-----------------------------------------------------------------
%---------------------------------Section 2-----------------------
%-----------------------------------------------------------------
\vskip 5mm

\section
{ $L_{-\frac{n}{n+2}}$ affine surface area of the polar body}

It was proved in \cite{MW2} that for a convex body $K\in C^2_+$

\begin{eqnarray}\label{as-1}
\lim _{\delta \rightarrow 0} \ c_n \frac{|(K_{\delta }
)^\circ|-|K^\circ|}{\delta ^{\frac{2}{n+1}}} &=& \int _{ S^{n-1}}
\frac{d\sigma(u)}{f_K(u)^\frac{1}{n+1} \ h_K(u)^{n+1}}
=\int_{\partial K} \frac {\kappa_{K}(x)^{\frac{n+2}{n+1}}} {\left
\langle x,N_{K}\left(x\right)\right\rangle ^{n+1}} d\mu_{K} (x)
\nonumber \\ \nonumber \\ &=&  as_{-n(n+2)}(K),
\end{eqnarray}
where  $c_n=2\left(\frac{|B_2^{n-1}|}{n+1}\right)^{\frac{2}{n+1}}$
and $K_{\delta}$ is the convex floating body  \cite{SW1}: The
intersection of all halfspaces $H^+$ whose defining hyperplanes
$H$ cut off a set of volume $\delta$ from $K$.

\vskip 3mm \noindent

%-------------------------Counterexample-------------------------
%----------------------------------------------------------------

Assumptions on the boundary of $K$ are needed in order that
(\ref{as-1})
holds.  \\
To see that, consider $B^n_{\infty}=\{ x \in
\mathbb{R}^n: \mbox{max}_{1 \leq i \leq n}|x_i| \leq 1 \}$.
As $\dis \kappa_{B^n_{\infty}}(x)=0$ a.e. on $\partial
B^n_{\infty}$,
$$\dis \int_{\partial \cn n} \frac
{\kappa_{B^n_{\infty}}(x)^{\frac{n+2}{n+1}} } {\left\langle
x,N_{B^n_{\infty}}\left(x\right)\right\rangle
^{n+1}}\,d\mu_{B^n_{\infty}} (x) =0.$$
However
\begin{equation}\label{ex1}
\dis \lim _{\delta \rightarrow 0}
\ c_n \frac{|((B^n_{\infty})_{\delta }
)^\circ|-|(B^n_{\infty})^\circ|}{\delta
^{\frac{2}{n+1}}}=\infty.
\end{equation}

Indeed, writing $K$ for $B_{\infty }^n$, we will construct a $0$-symmetric convex body
$K_1$ such that  $K_{\delta}\subseteq K_1 \subseteq K$. Then
$K^{\circ}\subseteq K_1^\circ \subseteq K_{\delta}^{\circ}$.
Therefore, to  show (\ref{ex1}), it is enough to show that
$$\dis \lim _{\delta \rightarrow 0}c_n
\frac{|K_1^\circ|-|K^\circ|}{\delta ^{\frac{2}{n+1}}}=\infty.$$
Let $R^{+}=\{(x_j)_{j=1}^n:  x_j\geq 0, 1 \leq j \leq n \}$. It is
enough to consider
$K^{+}=R^{+} \cap K$  and  to construct $(K_1)^+=K_1 \cap R^+$. \\
We define $(K_1)^+$ to be the intersection of $R^+$ with the
half-spaces $H_i^+$, $1 \leq i \leq n+1$, where
$H_i=\{(x_j)_{j=1}^n: x_i=1\}$,  $1 \leq i \leq n$,   and
$H_{n+1}=\left\{(x_j)_{j=1}^n: \sum _{j=1}^n
x_j=n-\left(n!\delta\right)^{\frac{1}{n}}\right\}$, $\delta >0$
sufficiently small. Notice that the hyperplane $H_{n+1}$
(orthogonal to the vector $(1, \dots, 1) $ )   cuts off a set of
volume exactly $\delta$ from $K$ and  therefore $K_\delta \subset
K_1$.
\par
Moreover, $\dis K_1^\circ$ can be written as  a  convex hull:
$$\dis K_1^\circ=\mbox{co}\bigg( \left\{ \pm e_i, 1\leq
i\leq n\right\}\cup \left\{ \frac{1}{s} \ (\varepsilon_1,
\dots,\varepsilon_n), \varepsilon_j=\pm 1, 1\leq j \leq n \right\}
\bigg),
$$
where $s=n-(n!\delta)^{\frac{1}{n}}$. Hence
$$|K_1^\circ|=\frac{2^n}{n!}\cdot
\frac{n}{n-(n!\delta)^{\frac{1}{n}}} $$
and therefore
\begin{eqnarray*}
\lim _{\delta \rightarrow
0}\frac{|K_1^\circ |-|K^\circ |}{\delta
^{\frac{2}{n+1}}}=    \frac{2^n}{n!} \ \lim _{\delta \rightarrow 0} \delta
^{\frac{-2}{n+1}}
\ \frac{(n!\delta)^{\frac{1}{n}}}{(n-(n!\delta)^{\frac{1}{n}})}
= \infty.
\end{eqnarray*}

\vskip 4mm

Now we  show
\vskip 3mm
\bt \label{as-2}
Let  $K$ be a convex body in $C_+^2$ such that  $0\in
\mbox{int}(K)$. Then
$$\dis \lim _{\delta \rightarrow 0} \ c_n
\frac{|(K_{\delta})^\circ|-|K^\circ|}{\delta
^{\frac{2}{n+1}}}= as_{-\frac{n}{n+2}}(K^\circ).
$$
\et
\vskip 3mm
\noindent
As a corollary of (\ref{as-1}) and Theorem \ref{as-2}  we get that for a convex body
$K\in C^2_+$
\begin{equation} \label{polar1}
as_{-n(n+2)}(K)=as_{-\frac{n}{n+2}}(K^\circ).
\end{equation}
This is a special case for $p=-n(n+2)$ of the formula
$as_p(K)=as_{\frac{n^2}{p}} (K^\circ)$ proved in \cite{Hug} for $p
> 0$. We will show in the next section that this formula holds
for all $p<0, p\neq -n$ for convex bodies with sufficiently smooth
boundary. For $p=0$ (and $K\in C^2_+$) the formula holds trivially
as $as_0(K)=n |K|$ and $as_{\infty}(K^\circ)=n |K|$ (see
\cite{SW5}).

\vskip 3mm For the proof of Theorem \ref{as-2} we need the
following lemmas.

%------------------------------Lemma1----------------------------------
%----------------------------------------------------------------------
\bl \label{Lemma-1} Let $K\in C^2_+$. Then for any
$x\in \pt K^\circ$, we have
\begin{eqnarray*}
\lim _{\d \rightarrow 0} \frac{\langle x,
N_{K^\circ}(x)\rangle }{n\
\d^{\frac{2}{n+1}}}\left[\left(\frac{\|x_{\d}\|}{\|x\|}\right)^n-1\right]
=\frac{\langle x, N_{K^\circ}(x)\rangle ^2}{c_n \
 (\k _{K^\circ} (x))^{\frac{1}{n+1}} }
 \end{eqnarray*}
 where
 $x_{\d}\in \pt (K_{\d})^{\circ}$ is  in the ray passing through $0$ and $x$.
 \el

\vskip 3mm
\noindent
{\bf Proof }
\par
\noindent Since $K$, and hence also $K_{\d}$, are in $C^2_+$ one
has that $K^\circ$ and $(K_{\d})^\circ$ are in $C^2_+$. Therefore,
for $x \in  \pt K^\circ$ there exists a unique $y\in \pt K$, such
that, $\langle x,y \rangle=1$, namely $\dis
y=\frac{N_{K^o}(x)}{\langle N_{K^o}(x), x \rangle}$. $y$ has outer
normal vector $\dis N_K(y)=\frac{x}{\|x\|}$ and
$\frac{1}{\|x\|}=\langle y, N_K(y) \rangle$.

\vskip 3mm Similarly, for $x_\delta \in \partial (K_{\d})^\circ$
there exists a unique $y_{\d}$ in $\partial K_\delta$  such that
$\langle x_{\d},y_{\d} \rangle=1$, namely $\dis
y_{\d}=\frac{N_{(K_{\d})^o}(x_{\d})}{\langle
N_{(K_{\d})^o}(x_{\d}), x_{\d} \rangle}$, $y_{\d}$ has outer
normal vector $\dis N_{K_{\d}}(y_{\d})=
\frac{x_{\d}}{\|x_{\d}\|}=\frac{x}{\|x\|}$ and
$\frac{1}{\|x_{\d}\|}=\langle y_{\d}, N_{K_{\d}}(y_{\d}) \rangle$.

\vskip 3mm Let $y' = [0,y]\cap \pt K_{\d}$ ($[z_1,z_2]$ denotes
the line segment from $z_1$ to $z_2$) and let  $y_{\d}'\in \pt K$
be such that $y_\delta = [0,y_{\d}'] \cap K_\delta$.

\vskip 3mm We have
\begin{eqnarray*}
\frac{1}{\|x\|}&=& \langle y, N_K(y) \rangle \geq \langle y_{\d}',
N_K(y) \rangle=\langle y_{\d}', \frac{x}{\|x\|} \rangle,\\
\frac{1}{\|x_{\d}\|}&=& \langle y_{\d}, N_{K_{\d}}(y_{\d}) \rangle
\geq \langle y', N_{{K_{\d}}}(y_{\d}) \rangle=\langle y',
\frac{x}{\|x\|} \rangle.
\end{eqnarray*}
Hence
\begin{eqnarray*}
\left[\left(\frac{\|x_{\d}\|}{\|x\|}\right)^n-1\right]&=&\left[\left(\frac{\langle
y, N_K(y) \rangle}{\langle y_{\d}, N_{K_{\d}}(y_{\d})
\rangle}\right)^n-1\right] \geq \left[\left(\frac{\langle y_{\d}',
\frac{x}{\|x\|} \rangle}{\langle y_{\d}, \frac{x}{\|x\|}
\rangle}\right)^n-1\right] \nonumber
=\left[\left(\frac{\|y_{\d}'\|}{\|y_{\d}\|}\right)^n-1\right],
\end{eqnarray*}
$$
\left[\left(\frac{\|x_{\d}\|}{\|x\|}\right)^n-1\right]=\left[\left(\frac{\langle
y, N_K(y) \rangle}{\langle y_{\d}, N_{K_{\d}}(y_{\d})
\rangle}\right)^n-1\right] \leq \left[\left(\frac{\langle y,
\frac{x}{\|x\|} \rangle}{\langle y', \frac{x}{\|x\|}
\rangle}\right)^n-1\right]
$$
\begin{equation}\label{upper}
=\left[\left(\frac{\|y\|}{\|y'\|}\right)^n-1\right]
\end{equation}
and therefore
$$\frac{\langle x, N_{K^\circ}(x)\rangle
}{n}\left[\left(\frac{\|y_{\d}'\|}{\|y_{\d}\|}\right)^n-1\right]\leq
\frac{\langle x, N_{K^\circ}(x)\rangle
}{n}\left[\left(\frac{\|x_{\d}\|}{\|x\|}\right)^n-1\right]\leq
\frac{\langle x, N_{K^\circ}(x)\rangle }{n}
\left[\left(\frac{\|y\|}{\|y'\|}\right)^n-1\right].$$ \vskip 3mm
We first consider the lower bound.
\begin{eqnarray*}
\lim _{\d \rightarrow 0} \frac{\langle x, N_{K^\circ}(x)\rangle
}{n\ \d
^{\frac{2}{n+1}}}\left[\left(\frac{\|x_{\d}\|}{\|x\|}\right)^n-1\right]
&\geq& \lim _{\d \rightarrow 0} \frac{\langle x,
N_{K^\circ}(x)\rangle }{\langle y_{\d}',
N_{K_{\d}}(y_{\d}')\rangle}\frac{\langle y_{\d}',
N_{K_{\d}}(y_{\d}')\rangle }{n\ \d
^{\frac{2}{n+1}}}\left[\left(\frac{\|y_{\d}'\|}{\|y_{\d}\|}\right)^n-1\right].
\end{eqnarray*}
As $\delta \rightarrow 0$, $y_\delta'\rightarrow y$. As $K$ is in $C_+^2$,
$N_{K_{\d}}(y_{\d}') \rightarrow N_{K}(y)$ as $\delta \rightarrow 0$.
\par
\noindent
Therefore $\lim _{\d \rightarrow 0} \langle y_{\d}',
N_{K_{\d}}(y_{\d}')\rangle= \langle y, N_{K}(y)\rangle$.
By Lemma 7 and Lemma 10 of \cite{SW1},
$$
\lim _{\d \rightarrow 0} \frac{\langle y_{\d}',
N_{K_{\d}}(y_{\d}')\rangle }{n\ \d
^{\frac{2}{n+1}}}\left[\left(\frac{\|y_{\d}'\|}{\|y_{\d}\|}\right)^n-1\right]
=\frac{(\kappa_K(y))^\frac{1}{n+1}}{c_n}.
$$
Hence
\begin{eqnarray*}
&&\lim _{\d \rightarrow 0} \frac{\langle x, N_{K^\circ}(x)\rangle
}{n\ \d
^{\frac{2}{n+1}}}\left[\left(\frac{\|x_{\d}\|}{\|x\|}\right)^n-1\right]
\geq \frac{\langle x, N_{K^\circ}(x)\rangle }{\langle y,
N_{K}(y)\rangle}
 \frac{(\k _K (y))^{\frac{1}{n+1}}}{c_n }
= \frac{\langle x, N_{K^\circ}(x)\rangle ^2}{c_n\
 (\k _{K^\circ} (x))^{\frac{1}{n+1}} }.
 \end{eqnarray*}
The last equation follows from the fact that if $K\in C^2_+$,
then, for any $y\in \pt K$, there is  a unique point $x\in \pt
K^\circ$ such that $\langle x,y \rangle =1$ and \cite{Hug}
\begin{equation} \label{curvatureformula}
 \dis \langle y, N_K(y) \rangle \langle x,
N_{K^\circ} (x) \rangle =\left(\k _K(y)\k
_{K^\circ}(x)\right)^{\frac{1}{n+1}}.
\end{equation}

\vskip 3mm
Similarly, one gets for  the upper bound
\begin{eqnarray*}&&
\lim _{\d \rightarrow 0} \frac{\langle x,
N_{K^\circ}(x)\rangle }{n\
\d^{\frac{2}{n+1}}}\left[\left(\frac{\|x_{\d}\|}{\|x\|}\right)^n-1\right]
\leq  \frac{\langle x, N_{K^\circ}(x)\rangle ^2}{c_n \
 (\k _{K^\circ} (x))^{\frac{1}{n+1}} },
 \end{eqnarray*}
\noindent
hence altogether
\begin{eqnarray*}\lim _{\d \rightarrow 0}
\frac{\langle x, N_{K^\circ}(x)\rangle }{n\
\d^{\frac{2}{n+1}}}\left[\left(\frac{\|x_{\d}\|}{\|x\|}\right)^n-1\right]
=\frac{\langle x, N_{K^\circ}(x)\rangle ^2}{c_n\
 (\k _{K^\circ} (x))^{\frac{1}{n+1}} }.
 \end{eqnarray*}

%-------------------Lemma 2------------------------
%--------------------------------------------------

\vskip 5mm
\bl \label{Lemma-2}
Let $K\in C^2_+$. Then we have
\begin{eqnarray*}
\frac{\langle x, N_{K^\circ}(x)\rangle }{n\
\d^{\frac{2}{n+1}}}\left[\left(\frac{\|x_{\d}\|}{\|x\|}\right)^n-1\right]
\leq c(K,n),
\end{eqnarray*}
where $c(K,n)$ is a constant (depending on $K$ and $n$ only) and
$x$ and $x_{\d}$ are as in Lemma \ref{Lemma-1}. \el \vskip 3mm
\noindent {\bf Proof }  By (\ref{upper})
\begin{eqnarray*}\frac{\langle x, N_{K^\circ}(x)\rangle }{n\ \d
^{\frac{2}{n+1}}}\left[\left(\frac{\|x_{\d}\|}{\|x\|}\right)^n-1\right]
&\leq & \frac{\langle x, N_{K^\circ}(x)\rangle }{\langle y,
N_{K}(y)\rangle}\frac{\langle y, N_{K}(y)\rangle }{n\ \d
^{\frac{2}{n+1}}}\left[\left(\frac{\|y\|}{\|y'\|}\right)^n-1\right]\\&\leq
&\frac{\langle x, N_{K^\circ}(x)\rangle }{\langle y,
N_{K}(y)\rangle}\left(\frac{\|y\|}{\|y'\|}\right)^n \frac{\langle
y, N_{K}(y)\rangle }{n\ \d ^{\frac{2}{n+1}}} \left[1-
\left(\frac{\|y'\|}{\|y\|}\right)^n\right].
\end{eqnarray*}
Since $K_{\d}$ is increasing to $K$ as $\delta \rightarrow 0$,
there exists $\d_0>0$ such that for all  $\d<\d_0$,  $0\in
\mbox{int}(K_{\d})$. Therefore there exits $\a>0$ such that
$B^n_2(0,\a)\subset K_{\d}\subset K\subset B^n_2(0, \frac{1}{\a})$ for all
$\d <\d _0$. $B^n_2(0, r)$ is the $n$-dimensional Euclidean ball centered at $0$ with radius $r$.

Hence for $\d <\d _0 $
\begin{eqnarray*}\frac{\langle x, N_{K^\circ}(x)\rangle }{n\ \d
^{\frac{2}{n+1}}}\left[\left(\frac{\|x_{\d}\|}{\|x\|}\right)^n-1\right]
&\leq & \a ^{-2(n+1)}\ \frac{\langle y, N_{K}(y)\rangle }{n\ \d
^{\frac{2}{n+1}}}\left[1-
\left(\frac{\|y'\|}{\|y\|}\right)^n\right] \leq C' r(y)
^{-\frac{n-1}{n+1}}
\end{eqnarray*}
due to Lemma 6 in \cite{SW1}.  Here $r(y)$ is the radius of the biggest Euclidean ball
contained in $K$ and touching $\partial K$ at $y$.

Since $K$ is $ C^2_+$, by the Blaschke rolling theorem (see
\cite{Sch})  there is $r_0 >0$ such that $\dis r_0 \leq \min
_{y\in \pt K }r(y)$. We put $c(K,n)=C'r_0^{-\frac{n-1}{n+1}}$.

\vskip 5mm

%--------------------------------------------------------------
%-----------------Proof of theorem 2.1-------------------------
%--------------------------------------------------------------
\noindent
{\bf Proof of Theorem \ref{as-2}.}

$$
\frac{|(K_{\d})^\circ|-|K^\circ|}{\d ^{\frac{2}{n+1}}}=\frac{1}{n \ \d ^{\frac{2}{n+1}}}\int _{\partial K^o}
\langle x, N_{K^\circ}(x)\rangle
\left[\left(\frac{\|x_{\d}\|}{\|x\|}\right)^n-1\right]\,d\mu_{K^\circ} (x).
$$
Combining Lemma \ref{Lemma-1}, Lemma \ref{Lemma-2} and Lebesgue's
convergence theorem, gives Theorem \ref{as-2}:
 \begin{eqnarray*}
 \lim _{\d
\rightarrow 0}\frac{|(K_{\d})^\circ|-|K^\circ|}{\d
^{\frac{2}{n+1}}}&=& \lim _{\d \rightarrow 0}\frac{1}{n \ \d
^{\frac{2}{n+1}}}\int _{\partial K^o} \langle x,
N_{K^\circ}(x)\rangle
\left[\left(\frac{\|x_{\d}\|}{\|x\|}\right)^n-1\right]\,d\mu_{K^\circ} (x)\\
&=& \int _{\partial K^o}  \lim _{\d \rightarrow 0}\frac{1}{n \ \d
^{\frac{2}{n+1}}}\langle x, N_{K^\circ}(x)\rangle
\left[\left(\frac{\|x_{\d}\|}{\|x\|}\right)^n-1\right]\,d\mu_{K^\circ}
(x)\\&=& \int _{\partial K^o} \frac{\langle x,
N_{K^\circ}(x)\rangle ^2}{c_n \
 (\k _{K^\circ} (x))^{\frac{1}{n+1}}  } \,d\mu_{K^\circ} (x)\\&=&\frac{1}{c_n}\  as_{\frac{-n}{n+2}}(K^\circ).
 \end{eqnarray*}
\vskip 5mm
%--------------------------------------------------------------
%----------------------------Remarks-------------------------
%--------------------------------------------------------------
\noindent
{\bf Remark}

The proof of Theorem \ref{as-2} provides a uniform method to
evaluate
$$\lim _{t \rightarrow 0}\frac{|(K_{t})^\circ|-|K^\circ|}{t
^{\frac{2}{n+1}}}$$ where $K_{t}$ is a family convex bodies
constructed from the convex body $K$ such that $K_{t}\subset K$ or- similarly-
such that  $K\subset
K_{t}$.
In particular,
we can apply this method  to prove the analog statements as in
(\ref{as-1}) and Theorem \ref{as-2}
 if we take as $K_t$ the illumination body  of $K$ \cite{W1}, or the Santal\'o body of $K$  \cite{MW1},
or the convolution body of $K$  \cite{ Schm} - and there are many more.
\vskip 5mm
%--------------------------------------------------------------
%----------------------------Section 3-------------------------
%--------------------------------------------------------------

\section{$L_p$ affine surface areas}

We now prove that for all $p \neq -n$ and all $K \in C^2_+$,
$as_{p}(K)=as_{\frac{n^2}{p}}(K^\circ)$.
To do so, we use the surface body of a convex body
which was introduced in \cite{SW4, SW5}.  We also give a new geometric interpretation of $L_p$
affine surface area for all $p \neq -n$.

\begin{definition}
Let $s\geq 0$ and $f:\partial K\rightarrow \mathbb{R}$ be a
nonnegative, integrable function. The {\it surface body $K_{f,s}$}
is the intersection of all the closed half-spaces $H^+$ whose
defining hyperplanes $H$ cut off a set of $f \mu_K$-measure less
than or equal to $s$ from $\partial K$. More precisely,
\begin{equation*}\label{def1}
K_{f,s}= \bigcap _{\int_{\partial K \cap H^-}f d\mu_K \leq s} H^+.
\end{equation*}
\end{definition}
\vskip 3mm
\bt \label{asp}
Let  $K$ be a convex body in $C^2_+$
and such that  $0$ is the center of gravity of $K$. Let
$f:\partial K\rightarrow \bbR$ be an integrable function such that
$f(x)>c$ for all $x\in \pt K$ and some constant $c>0$. Let
$\beta_n=2\left(|B^{n-1}_2 |\right)^{\frac{2}{n-1}}$. Then
\begin{eqnarray*}
\lim _{s\rightarrow 0} \beta_n
\frac{|(K_{f,s})^\circ|-|K^\circ|}{s^{\frac{2}{n-1}}}=\int_{S^{n-1}}\frac{d\sigma(u)}{h_K(u)^{n+1}
f_K(u)^\frac{1}{n-1} \big(f(N_K^{-1}(u))\big)^\frac{2}{n-1}}
\end{eqnarray*}
where $N_K: \partial K \rightarrow S^{n-1}$, $x \rightarrow N_K(x)=u$ is
the Gauss map.

\et
\vskip 3mm
\noindent
{\bf Proof}
\par
\noindent Let $u\in S^{n-1}$.  Let $x \in \partial K$ be such that
$N_K(x)=u$ and let $x_s \in \partial K_{f,s}$ be such that
$N_{K_{f,s}}(x_s)=u$. Let $H_{\Delta}=H(x-\Delta u,u)$ be the
hyperplane  through $x-\Delta u$ with outer normal vector $u$.
Since $K$ has everywhere strictly positive Gaussian curvature, by
Lemma 21 in \cite{SW5} almost everywhere on $\partial K$,
$$\lim _{\Delta \rightarrow 0}
\frac{1}{|\partial K\cap H_{\Delta}^-|}\int _{\partial K\cap
H_{\Delta}^-} |f(x)-f(y)|\,d\mu_{K}(y)=0.
$$ This implies  that
\begin{equation}\label{f-lim}
\lim _{\Delta \rightarrow 0} \frac{1}{|\partial K\cap
H_{\Delta}^-|}\int _{\partial K\cap H_{\Delta}^-}
f(y)\,d\mu_{\partial K}(y)=f(x).
\end{equation}
Let $b_s=h_{K}(u)-h_{K_{f,s}}(u)$. As  $H(x-b_s u, u)
=H(x_s,u)$ (the hyperplane through $x_s$ with outer normal
$u$)  and as $b_s \rightarrow 0$ as $s \rightarrow 0$,  (\ref{f-lim}) implies
\begin{equation}\label{f-lima}
\lim _{s\rightarrow 0} \frac{1}{|\partial K\cap H^-(x_s,u)|}\int
_{\partial K\cap H^-(x_s,u)} f(y)\,d\mu_{ K}(y)=f(x).
\end{equation}
Hence there exists   $s_1$ small
enough, such that for all $s<s_1$,
\begin{equation}\label{eq-1}
s  \leq \int _{\partial K\cap H^-(x_s,u)} f(y)\,d\mu_{ K}(y) \leq
(1+\varepsilon ) f(x){|\partial K\cap H^-(x_s,u)|}.
\end{equation}

As  $\partial K$ has everywhere strictly  positive Gaussian
curvature, the indicatrix of Dupin exists everywhere on $\partial
K$ and  is an ellipsoid. It then follows from (\ref{eq-1}) with
Lemmas 1.2, 1.3 and 1.4  in \cite{SW4} that there exists $0< s_2 <
s_1$ such that for all $0<s<s_2$

$$
s\leq (1+ \varepsilon )\ f(x)\  |B^{n-1}_{2}| \ {\sqrt{f _K(u)}}\
\left(2 b_s \right)^{\frac{n-1}{2}},
$$
or, equivalently
\begin{equation}\label{eq-5}
\frac{b_s}{s^{\frac{2}{n-1}}} \geq \frac{1- c_1 \
\varepsilon }{\beta _n \ f(N_K^{-1}(u))^{\frac{2}{n-1}} \ f_K
(u)^{\frac{1}{n-1}}},
\end{equation}
where $c_1$ is an absolute constant.

%--------------------------------------------------------------------------------------
%-------------------From above---------------------------------------------------------
%-------------------------------------------------------------------------------------
\vskip 3mm
Let now $x_s' \in [0, x]\cap \pt
K_{f,s}$. Then $\langle x_s', u \rangle \leq
h_{K_{f,s}}(u)$. Therefore $b_s=h_K(u)-h_{K_{f,s}}(u)\leq \langle
x-x_s', u\rangle$.

Hence for $s$ sufficiently small
\begin{eqnarray}\label{eq-5a}
\frac{b_s}{s^{\frac{2}{n-1}}} \leq \frac{\langle x-x_s', u\rangle
}{s^{\frac{2}{n-1}}} &\leq& \frac{\langle x,
u\rangle}{s^{\frac{2}{n-1}}} \left( 1-\frac{\| x_s'\|}{\|
x\|}\right) \leq
 \frac{\langle x, u\rangle}{s^{\frac{2}{n-1}}}
\  \frac{\|x_s' -x\|}{\|x\|}   \nonumber
\\\nonumber \\ &\leq& (1+\varepsilon)\
\frac{\langle x, N_K(x)
\rangle }{n \ s^{\frac{2}{n-1}} } \left[1-\left(\frac{\|x_s'\|}{\|x\|}\right)^n\right].
\end{eqnarray}
The last inequality follows as
$
1-\left(\frac{\|x_s'\|}{\|x\|}\right)^n \geq (1-\varepsilon) \  \frac{n\  \|x_s' -x\|}{\|x\|}
$ for  sufficiently small $s$.
By Lemma 23 in \cite{SW5}
\begin{equation}\label{eq-6}
\lim _{s\rightarrow 0 }\frac{1}{n s^{\frac{2}{n-1}} }\langle x,
N_K(x) \rangle
\left[1-\left(\frac{\|x_s'\|}{\|x\|}\right)^n\right]
 = \frac{1}{\beta _n
f(N_K^{-1}(u))^{\frac{2}{n-1}}f_K (u)^{\frac{1}{n-1}}}.
\end{equation}
Thus we get  from (\ref{eq-5}), (\ref{eq-5a}) and (\ref{eq-6}) that
\begin{equation}\label{eq-6a}
\lim _{s\rightarrow 0}\frac{b_s}{s^{\frac{2}{n-1}}}=  \frac{1}{\beta _n
f(N_K^{-1}(u))^{\frac{2}{n-1}}f_K (u)^{\frac{1}{n-1}}}.
\end{equation}

As $(1-t)^{-n}\geq 1+nt$ for  $0 \leq t <1$ and by (\ref{eq-6a}),
\begin{eqnarray}\label{eq-9}
\lim _{s\rightarrow 0}\frac{ \beta _n }{ns^{\frac{2}{n-1}}}\left(
[h_{K_{f,s}}(u)]^{-n} - [h_K(u)]^{-n} \right)&=&\lim
_{s\rightarrow 0}\frac{ \beta _n }{ns^{\frac{2}{n-1}}}
[h_K(u)]^{-n} \left[ \left(1+\frac{b_s}{h_K(u)}\right)^{-n}-1
\right]\nonumber \\ &\geq& \lim _{s\rightarrow 0}\frac{ \beta _n
}{[h_K(u)]^{n+1}} \frac{b_s}{s^{\frac{2}{n-1}}}\nonumber \\
&=& \frac{1}{[h_{K}(u)]^{n+1} \ f(N_K^{-1}(u))^{\frac{2}{n-1}}\ f_K
(u)^{\frac{1}{n-1}}}.
\end{eqnarray}
\par
\noindent
As $h_{K_{f,s}}(u)\geq \langle x_s', u\rangle$,
\begin{eqnarray}\label{eq-9a}
\frac{h_{K_{f,s}}(u)}{h_K(u)}\geq \frac{\langle x_s', u\rangle}{\langle x,u \rangle}=\frac{\|x_s'\|}{\|x\|}.
\end{eqnarray}
Since $K\in C^2_+$,
$h_{K_{f,s}}(u)\rightarrow h_K(u)$ as $s\rightarrow 0$.
Therefore,
$$
\lim _{s\rightarrow 0}\frac{ \beta _n
}{ns^{\frac{2}{n-1}}}\left([ {h_{K_{f,s}}(u)}]^{-n} -
[h_K(u)]^{-n} \right)=\lim _{s\rightarrow 0}\frac{ \beta _n
}{ns^{\frac{2}{n-1}}}[{h_{K_{f,s}}(u)}]^{-n}\left[1 -\left(
\frac{h_{K_{f,s}}(u)}{h_K(u)}\right)^{n}\right]
$$
$$
 \leq \lim _{s\rightarrow 0}\frac{ \beta _n
}{ns^{\frac{2}{n-1}}}[{h_{K_{f,s}}(u)}]^{-n}\left[1 -\left(
\frac{\|x_s'\|}{\|x\|}\right)^{n}\right]
= \lim _{s\rightarrow 0}\frac{ \beta _n
}{ns^{\frac{2}{n-1}}}[{h_{K_{f,s}}(u)}]^{-n}\frac{\langle x,u
\rangle }{\langle x,u \rangle}\left[1 -\left(
\frac{\|x_s'\|}{\|x\|}\right)^{n}\right]
$$
$$
=\lim _{s\rightarrow 0} \frac{1}{h_K(u) \ [{h_{K_{f,s}}(u)}]^{n}}
\ \lim _{s\rightarrow 0}\frac{ \beta _n
}{ns^{\frac{2}{n-1}}}\langle x,N_K(x) \rangle\left[1 -\left(
\frac{\|x_s'\|}{\|x\|}\right)^{n}\right]
$$
\begin{equation}\label{eq-8}
=\frac{1}{[h_K(u)]^{n+1} f(N_K^{-1}(u))^{\frac{2}{n-1}}f_K
(u)^{\frac{1}{n-1}}}
\end{equation}
where the last equality follows from (\ref{eq-6}). \\
\par
\noindent
Altogether, (\ref{eq-9})  and (\ref{eq-8}) give
\begin{equation}\label{eq-12}
\lim _{s\rightarrow 0}\frac{ \beta _n }{ns^{\frac{2}{n-1}}}\left([
{h_{K_{f,s}}(u)}]^{-n} - [h_K(u)]^{-n} \right)\nonumber =
\frac{1}{[h_{K}(u)]^{n+1} \ f(N_K^{-1}(u))^{\frac{2}{n-1}}\ f_K
(u)^{\frac{1}{n-1}}}.
\end{equation}
\par
\noindent
Therefore
\begin{eqnarray*}
\lim _{s\rightarrow 0} \beta_n
\frac{|(K_{f,s})^\circ|-|K^\circ|}{s^{\frac{2}{n-1}}} &=& \lim
_{s\rightarrow 0} \frac{\beta_n }{n \ s^{\frac{2}{n-1}}} \int
_{S^{n-1}}
\left[\left(\frac{1}{h_{K_{f,s}}(u)}\right)^n-\left(\frac{1}{h_{K}(u)}\right)^n\right]\,d\sigma(u)\\
&=& \int _{S^{n-1}} \lim _{s\rightarrow 0} \frac{\beta_n }{n \
s^{\frac{2}{n-1}}}
\left[\left(\frac{1}{h_{K_{f,s}}(u)}\right)^n-\left(\frac{1}{h_{K}(u)}\right)^n\right]\,d\sigma(u)\\
&=& \int_{S^{n-1}}\frac{d\sigma(u)}{h_K(u)^{n+1}
f_K(u)^\frac{1}{n-1} \big(f(N_K^{-1}(u))\big)^\frac{2}{n-1}},
\end{eqnarray*}
provided we can interchange integration and limit.
\vskip 3mm

We show this next.
%--------------------------------------------------------------
%----------------------------Interchanging order --------------
%--------------------------------------------------------------
To do
so, we show that for all $u \in S^{n-1} $ and all sufficiently small $s>0$,
$$
\frac{1}{n \ s^\frac{2}{n-1}}\  \bigg[
\bigg(\frac{1}{h_{K_{f,s}}(u)}\bigg)^n - \bigg(\frac{1}{h_{K}(u)}
\bigg)^n \bigg] \leq g(u)
$$
with $\int_{S^{n-1}} g(u) \ d\sigma (u)  < \infty$. As $0\in
int(K)$, the interior of $K$, there exists  $\a >0$ such that for
all $s$ sufficiently small $B^n_2 (0,\alpha) \subset K_{f,s}
\subset K \subset B^n_2(0,\frac{1}{\alpha})$. Therefore, $\a \leq
h_{K_{f,s}}(u)\leq h_K(u) \leq \frac{1}{\a}$ and $\a \leq
\frac{1}{h_{K}(u)}\leq \frac{1}{h_{K_{f,s}}(u)}\leq \frac{1}{\a}$.
\par
With (\ref{eq-9a}), we thus get for all $s>0$,
\begin{eqnarray*}
&&\frac{1}{n \ s^\frac{2}{n-1}}\  \bigg(
\big({h_{K_{f,s}}(u)}\big)^{-n} - \big({h_{K}(u)} \big)^{-n} \bigg)\\
& &=\frac{1}{n \ s^\frac{2}{n-1}} \big({h_{K_{f,s}}(u)}\big)^{-n}
\bigg(1-\frac{\big(h_{K_{f,s}}(u)\big)^n}{\big(h_{K}(u)\big)^n}
\bigg)\\& &\leq \frac{\a ^{-n}}{n \ s^\frac{2}{n-1}}  \bigg[1-
\bigg(\frac{\|x'_s\|}{\|x\|}\bigg)^n \bigg]
\\ & &\leq \a
^{-(n+1)}\ \ \frac{\langle x,u\rangle}{n \ s^\frac{2}{n-1}}
\bigg[1- \bigg(\frac{\|x'_s\|}{\|x\|}\bigg)^n  \bigg].
\end{eqnarray*}
\par
\noindent
By Lemma 17 in \cite{SW5} there exists $s_3 $ such that  for all $s \leq s_3$
\begin{eqnarray*}
\frac{\langle x,u\rangle}{ s^\frac{2}{n-1}}
\bigg[1- \bigg(\frac{\|x'_s\|}{\|x\|}\bigg)^n  \bigg] \leq
\frac{C}{(M_f(x))^{\frac{2}{n-1}}r(x)},
\end{eqnarray*}
where $C$ is an absolute constant and,  as in the proof of Lemma \ref{Lemma-2},  $r(x)$
is the biggest Euclidean ball contained in $K$ that  touches
$\partial K$ at $x$. Thus, as $\partial K$ is $C^2_+$, by Blaschke's rolling theorem (see \cite{Sch})
there is $r_0$ such that $r(x) \geq r_0$.
\par
$$M_f(x) = \inf _{0 <s}
\frac{\int_{\partial K \cap H^-(x_s, N_{K_{f,s}}(x_s))}
 f d\mu_K}{|\partial K \cap H^-(x_s, N_{K_{f,s}}(x_s)|}$$  is the {\em minimal function}. It was introduced in
 \cite{SW5}. By the assumption on $f$, $M_f(x) \geq c$. Thus altogether
\begin{eqnarray*}
\frac{1}{n \ s^\frac{2}{n-1}}\  \bigg(
\big({h_{K_{f,s}}(u)}\big)^{-n} - \big({h_{K}(u)} \big)^{-n}
\bigg) \leq \frac{\alpha ^{-(n+1)} \ C}{n\ c^\frac{2}{n-1} \
r_0}=g(u),
\end{eqnarray*}
which, as a constant, is integrable.

\vskip 5mm
 \bt \label{asp1}
 Let  $K$ be a convex body in
$C^2_+$ and such that  $0$ is the center of gravity of $K$. Let
$f:\partial K\rightarrow \bbR$ be an integrable function such that
$f(y)>c$ for all $y\in \pt K$ and some constant $c>0$. Let
$\beta_n=2\left(|B^{n-1}_2 |\right)^{\frac{2}{n-1}}$. Then
\begin{eqnarray*} \lim _{s\rightarrow 0} \beta_n
\frac{|(K_{f,s})^\circ|-|K^\circ|}{s^{\frac{2}{n-1}}}=\int_{\pt
K^\circ } \left(\frac{\langle x, N_{K^\circ }(x)\rangle}{\langle
y(x), N_K(y(x)) \rangle } \right) \left(\frac{\k
_K(y(x))^\frac{1}{n-1}}{f(y(x))^\frac{2}{n-1}}\right)\,d\mu_{K^\circ}(x)
\end{eqnarray*}
Here $y(x) \in  \pt K$ is such that $\langle y(x), x
\rangle=1 $.
\et
{\bf Proof}

We follow the pattern of the proof of Theorem \ref{asp}
integrating now  over $\pt K^\circ$ instead of $S^{n-1}$.

\vskip 5mm

As a corollary we get the following geometric interpretation
of $L_p$ affine surface area.

\vskip 3mm
\begin{corollary} \label{general-p-affine-surface1}
 Let $K \in C^2_+$ be a convex body. For $p \in \mathbb{R}$,  $p\neq -n$, let
 $f_p: \partial K \rightarrow {\mathbb R}$ be defined by $f_p(y)=\k_K(y)^\frac{n^2+p}{2(n+p)}\
  \langle y, N_K(y)\rangle^\frac{-(n-1)(n^2+2n+p)}{2(n+p)}$. Then
\par
(i) $$ \lim _{s\rightarrow 0}\beta _n \ \frac{|(K
_{f_p,s})^\circ|-|K^\circ|}{s^{\frac{2}{n-1}}}=
as_{\frac{n^2}{p}}(K^\circ).
$$
\par
(ii) $$ \lim _{s\rightarrow 0}\beta _n \ \frac{|(K
_{f_p,s})^\circ|-|K^\circ|}{ s^{\frac{2}{n-1}}}= as_{p}(K).
$$
\par
(iii) $$ as_{p}(K)=as_{\frac{n^2}{p}}(K^\circ).
$$
\end{corollary}

\vskip 3mm
\noindent
{\bf Proof}
\par
\noindent
Notice first that $f_p(y)$ verifies the conditions of  Theorems  \ref{asp} and \ref{asp1}.

(i) For $x \in \partial K^\circ$, let now $y(x)$ be the (unique) element in $\partial K$ such that $\langle x,y(x) \rangle =1$.
Then,  by Theorem  \ref{asp1}, with $f(y(x))=f_p(y(x))$, and with (\ref{curvatureformula})

\begin{eqnarray*} \label{gl1}
\lim _{s\rightarrow 0} \beta _n
\frac{|\left(K_{f_p,s}\right)^\circ|-|K^\circ|}{s^{\frac{2}{n-1}}}&=&
  \int _{\pt K ^\circ} \langle x, N_{K^\circ}(x)\rangle
 \frac{\langle y(x), N_K(y(x))\rangle^\frac{n(n+1)}{n+p}}{\kappa_K(y(x))^\frac{n}{n+p}} \,d\mu _{K^\circ}
 (x)\\ \\
 &=&
  \int _{\partial K^\circ} \frac{\kappa_{K^\circ}(x)^{\frac{n}{n+p}} }{\langle x, N_{K^\circ}(x)\rangle
^{\frac{n^2-p}{n+p}}} \,d\mu _{K^\circ} (x)
=as_\frac{n^2}{p}(K^\circ).
\end{eqnarray*}

\par
(ii) For $u \in S^{n-1}$, let now $y \in \partial K$ be such that
$N_K(y)=u$. Then $f_p(N_K^{-1}(u))=f_K(u)^{-\frac{n^2+p}{2(n+p)}}\
h_K(u)^\frac{-(n-1)(n^2+2n+p)}{2(n+p)}$. By Theorem   \ref{asp}
with $f(N_K^{-1}(u))=f_p(N_K^{-1}(u))$

\begin{eqnarray*} \label{gl2}
\lim _{s\rightarrow 0} \beta
_n\frac{|\left(K_{f_p,s}\right)^\circ|-|K^\circ|}{s^{\frac{2}{n-1}}}
=\int _{S^{n-1}}\frac{f _K(u)^{\frac{n}{n+p}}
}{h_K(u)^{\frac{n(p-1)}{n+p}}}\,d\sigma (u)=as_{p}(K).
\end{eqnarray*}
\par
(iii)  follows from (i) and (ii).

\vskip 3mm

%--------------------------------------------------------------
%----------------------------Section 4-------------------------
%--------------------------------------------------------------
\section{Inequalities}

\vskip 3mm
\noindent

%-------------------------------Theorem 4.1------------------------
%------------------------------------------------------------------
\bt\label{p-aff-sio8}
Let $s\neq -n,  r \neq -n, t \neq -n$ be
real numbers. Let $K$ be a convex body in $\mathbb R^n$ with
centroid  at the origin and such that $\mu _K \{x \in \partial K:
\kappa_K(x) =0\}=0$.
\par
(i) If $\frac{(n+r)(t-s)}{(n+t)(r-s)}>1$, then
\begin{equation}\label{i-2} \nonumber
as_r(K)\leq \big(as_t(K) \big)^{\frac{(r-s)(n+t)}{(t-s)(n+r)}}
\big(as_s (K)\big)^{\frac{(t-r)(n+s)}{(t-s)(n+r)}}.
\end{equation}

\par
(ii) If $\frac{(n+r)t}{(n+t)r} >1$,  then
$$
\bigg(\frac{as_r(K)}{n|K|}\bigg)\leq
\bigg(\frac{as_t(K)}{n|K|}\bigg)^{\frac{r(n+t)}{{t(n+r)}}}.
$$
\et

\par
\vskip 2mm \noindent {\bf Proof}
\par
(i)  By
H\"older's inequality -which enforces the condition $\frac{(n+r)(s-t)}{(n+t)(s-r)} >1$
 \begin{eqnarray} \label{firstcase}
 as_r(K)&=&\int _{\pt
K} \frac{\k _K(x)^{\frac{r}{n+r}}}{\langle x, N_K(x)\rangle
^{\frac{n(r-1)}{n+r}}}\, d\mu _K(x) \nonumber
\\
&=& \int _{\pt K}
\left(\frac{\k _K(x)^{\frac{t}{n+t}}}{\langle x, N_K(x)\rangle
^{\frac{n(t-1)}{n+t}}} \right)^{\frac{(r-s)(n+t)}{(t-s)(n+r)}}
\left(\frac{\k _K(x)^{\frac{s}{n+s}}}{\langle x, N_K(x)\rangle
^{\frac{n(s-1)}{n+s}}}\right)^{\frac{(t-r)(n+s)}{(t-s)(n+r)}}\,
d\mu _K(x) \nonumber
\\ &\leq & \big(as_t(K)
\big)^{\frac{(r-s)(n+t)}{(t-s)(n+r)}} \big(as_s
(K)\big)^{\frac{(t-r)(n+s)}{(t-s)(n+r)}}. \nonumber
 \end{eqnarray}
\par
(ii) Similarly, again using H\"older's inequality -which now
enforces the condition $~\frac{(n+r)t}{(n+t)r} >1,~$
\begin{eqnarray*}
 as_r(K)&=&\int _{\pt
K} \frac{\k _K(x)^{\frac{r}{n+r}}}{\langle x, N_K(x)\rangle
^{\frac{n(r-1)}{n+r}}}\, d\mu _K(x)  = \int _{\pt K}
\left(\frac{\k _K(x)^{\frac{t}{n+t}}}{\langle x, N_K(x)\rangle
^{\frac{n(t-1)}{n+t}}} \right)^{\frac{r(n+t)}{t(n+r)}} \frac{d\mu
_K(x)}{\langle x, N_K(x)\rangle ^{\frac{(r-t)n}{(n+r)t}}} \\ \\
&\leq &\big(as_t(K)\big)^{\frac{r(n+t)}{t(n+r)}} \big(n\
|K|\big)^{\frac{(t-r)n}{(n+r)t}}.
\end{eqnarray*}

\vskip 3mm Condition $\frac{(n+r)(t-s)}{(n+t)(r-s)}>1$ implies $8$
cases: $-n<s<r<t$, $s<-n<t<r$,  $r<t<-n<s$, $t<r<s<-n$,
$s<r<t<-n$,  $r<s<-n<t$, $t<-n<s<r$ and $-n<t<r<s$.

\vskip 3mm Note also that (ii) describes a monotonicity condition
for $\left(\frac{as_r(K)}{n|K|}\right)^{\frac{{n+r}}{r}}$: if $0 <
r <t$, or $r<t<-n$, or $-n<r<t<0$ then
$$
\left(\frac{as_r(K)}{n|K|}\right)^{\frac{{n+r}}{r}} \leq
\left(\frac{as_t(K)}{n|K|}\right)^{\frac{{n+t}}{t}}.
$$
\vskip 3mm We now analyze various subcases of Theorem
\ref{p-aff-sio8} (i) and (ii). For $r=0$, if
$\frac{n(s-t)}{s(n+t)}
>1$

\begin{equation*}  \label{4.1.1}
n|K| \leq \big(as_t(K) \big)^{\frac{s(n+t)}{n(s-t)}} \big(as_s
(K)\big)^{\frac{t(n+s)}{n(t-s)}}.
\end{equation*}
For $s =0$, if
 $\frac{t(n+r)}{r(n+t)} >1$,
\begin{equation}  \label{4.1.2}
as_r(K) \leq
\left(n|K| \right)^\frac{n(t-r)}{t(n+r)} \big(as_t(K)
\big)^{\frac{r(n+t)}{t(n+r)}}.
\end{equation}
For $s \rightarrow \infty$, if $\frac{n+r}{n+t} >1$,
\begin{equation}  \label{4.1.3}
as_r(K) \leq \big(as_{\infty}(K)\big)^\frac{r-t}{n+r} \big(as_t(K)
\big)^{\frac{n+t}{n+r}}.
\end{equation}
For  $r\rightarrow \infty$, if
$\frac{t-s}{n+t}>1$ and if $K$ is in $C^2_+$,
\begin{equation} \label{4.1.4}
as_{\infty}(K)=n|K^\circ| \leq \big(as_t(K)
\big)^{\frac{n+t}{t-s}} \big(as_s (K)\big)^{\frac{n+s}{s-t}}.
\end{equation}

As for all convex bodies $K$, $as_{\infty}(K)\leq n |K^\circ|$
(see \cite{SW5}), it follows from (\ref{4.1.2}) that, for all
convex body $K$ with centroid at origin,
\begin{equation}  \label{4.1.5a}
as_r(K)\leq (n |K|)^{\frac{n}{n+r}}\  (n
|K^\circ|)^{\frac{r}{n+r}}, \hskip 3mm r>0
\end{equation}
and from (\ref{4.1.3}),
\begin{equation}  \label{4.1.5b}
n |K|  \big(n |K^\circ |\big)^\frac{t}{n}  \leq \big(as_t(K)
\big)^{\frac{n+t}{n}}, \hskip 3mm -n < t <0.
\end{equation}
Similarly, (\ref{4.1.4}) implies that,  if in addition $K$ is in $ C^2_+$,
\begin{equation}  \label{4.1.5c}
n |K^\circ| \big(n|K|\big)^{\frac{n}{t}} \leq  \big(as_t(K)
\big)^{\frac{n+t}{t}}, \hskip 3mm t < -n
\end{equation}
(\ref{4.1.5a}) can also be obtained from
Proposition 4.6 of \cite{Lu1} and Theorem 3.2 of \cite{Hug}.

%-------------------------------Corollary 4.1------------------------
%------------------------------------------------------------------
\vskip 5mm

Inequalities (\ref{4.1.5a}),  (\ref{4.1.5b}) and
(\ref{4.1.5c}) yield the following Corollary which was proved
by Lutwak \cite{Lu1} in the  case $p \geq 1$.

\begin{corollary} \label{p-aff-iso3}
Let $K$ be convex body in $\mathbb R^n$ with centroid at the
origin.

(i) For all $p\geq 0$
\begin{equation*}
as_p(K)\ as_p(K^\circ)\leq n^2|K|\ |K^\circ|.
\end{equation*}

(ii) For $-n<p<0$,
$$
as_p(K)as_p(K^\circ)\geq n^2|K|\ |K^\circ|.
$$
If $K$ is in addition in $C^2_+$,  inequality (ii) holds for all
$p<-n$.
\end{corollary}
\vskip 3mm Thus, using Santal\'o inequality in (i), for $ p\geq
0$, $ as_p(K)as_p(K^\circ)\leq  as_p(\bn n)^2$, and inverse
Santal\'o inequality in (ii), for $ -n<p<0$,
$as_p(K)as_p(K^\circ)\geq  c^n as_p(\bn n)^2$. $c$ is the constant
in the inverse Santal\'o inequality \cite{BM,GK2}.

\vskip 3mm
\noindent
{\bf Proof}
\par
\noindent
(i) follows immediately form (\ref{4.1.5a}).
(ii) follows from (\ref{4.1.5b}) if $-n<p<0$ and from
(\ref{4.1.5c}) if $p<-n$.

\vskip 5mm
\noindent

Lutwak  \cite{Lu1} proved   for $p \geq 1 $
 \begin{equation*}\label{Lut2}
\frac{as_p(K)}{as_p(B^n_2)}\leq  \left(\frac{|K|}{|\bn
n|}\right)^{\frac{n-p}{n+p}}
\end{equation*} with equality if and only if $K$ is an ellipsoid.
We now generalize these $L_p$-affine isoperimetric inequalities  to $p <1$.

%-------------------------------Theorem 4.2------------------------
%------------------------------------------------------------------

\vskip 5mm

\bt\label{p-aff-iso4}
\par
 Let  $K$ be a convex body with
centroid at the origin.

(i) If $p\geq 0$, then
\begin{eqnarray*}
\frac{as_p(K)}{as_p(B^n_2)}\leq \left(\frac{|K|}{|\bn
n|}\right)^{\frac{n-p}{n+p}},
\end{eqnarray*}
with equality if and only if $K$ is an ellipsoid. For $p=0$, equality holds trivially for all $K$.
\par
(ii) If
$-n<p<0$, then
\begin{eqnarray*}
\frac{as_p(K)}{as_p(B^n_2)}\geq
\left(\frac{|K|}{|B^n_2|}\right)^{\frac{n-p}{n+p}},
\end{eqnarray*} with
equality if and only if $K$ is an ellipsoid.
\par
(iii) If  $K$ is in addition  in $C^2_+$ and if $p < -n$, then
\begin{equation*}
\frac{as_p(K)}{as_p(B^n_2 )}\geq
c^{\frac{np}{n+p}}\left(\frac{|K|}{|\bn
n|}\right)^{\frac{n-p}{n+p}}.
\end{equation*}
where $c$ is the constant in the inverse Santal\'o inequality
\cite{BM,GK2}. \et

\vskip 3mm
We cannot expect to get a strictly positive lower bound in Theorem \ref{p-aff-iso4} (i),
 even if $K$ is in
$C^2_+$: Consider, in ${\mathbb R}^2$, the convex body $K(R, \varepsilon)$ obtained
as the intersection of four
Euclidean balls with radius $R$ centered at $(\pm(R-1),0)$, $(0, \pm(R-1))$, $R$ arbitrarily large.
To obtain a body in $C^2_+$, we ``round" the corners  by putting there Euclidean balls with radius
 $\varepsilon$, $\varepsilon$ arbitrarily small. Then  $as_p(K(R, \varepsilon)) \leq \frac{16}{R^\frac{p}{2+p}} + 4  \pi\ \varepsilon^\frac{2}{2+p}$. A similar construction can be done in higher  dimensions.
\par
This example also shows that, likewise,
we cannot expect  finite upper bounds in Theorem \ref{p-aff-iso4} (ii) and (iii).  If $-2 <p <0$, then  $as_p(K(R, \varepsilon))
 \geq 2 ^\frac{3(p+1)}{2+p}R^{\frac{-p}{2+p}}$. If $p < -2$, then $-2 < \frac{4}{p} <0$ and thus
$$
as_{p}(K(R,\varepsilon)^\circ)
 =as_{\frac{4}{p}}(K(R,\varepsilon))\geq R^{\frac{-2}{p+2}}\
 2^{\frac{12+3p}{4+2p}}.
 $$

\par
Note also that in part (iii)  we cannot remove the constant
$c^{\frac{np}{n+p}}$. Indeed, if $p\rightarrow -\infty$,  the
inequality
 becomes $c^n |B^n_2|^2 \leq |K| |K^\circ|$.

\vskip 3mm
\noindent
{\bf Proof of Theorem \ref{p-aff-iso4}}

(i) The case $p=0$ is  trivial. We prove the case $p>0$. Combining
inequality (\ref{4.1.5a}), the Blaschke Santal\'o inequality, and
$as_q(B^n_2)= n|B^n_2|^{\frac{n}{n+q}}|B^n_2|^{\frac{q}{n+q}}$,
one obtains
\begin{eqnarray*}
\frac{as_p(K)}{as_p(B^n_2)}\leq
\left(\frac{|K^\circ|}{|B^n_2|}\right)^{\frac{p}{n+p}}
\left(\frac{|K|}{|B^n_2|}\right)^{\frac{n}{n+p}} \leq
\left(\frac{|K|}{|B^n_2|}\right)^{\frac{n-p}{n+p}}.
\end{eqnarray*}
This proves the inequality. The equality case follows from the equality
case for the Blaschke Santal\'o inequality.

\vskip 2mm (ii) Combining inequality (\ref{4.1.5b}) and
$\big(as_p(\bn n)\big)^{\frac{n+p}{n}}=n |\bn n|  \big(n |\bn n
|\big)^\frac{p}{n},$ one gets, for $-n<p<0$,
\begin{eqnarray*}
\left(\frac{as_p(K)}{as_p(\bn n)}\right)^{\frac{n+p}{n}} \geq
\left(\frac{|K|}{|\bn n|}\right) \ \left(\frac{|K^\circ|}{|\bn n
|}\right)^\frac{p}{n} \geq \left(\frac{|K|}{|\bn n|}\right)
^\frac{n-p}{n}
\end{eqnarray*}
where the last inequality follows from the Blaschke Santal\'o
inequality. As $\frac{p}{n}<0$,
$$
(|K|\ |K^\circ|)^{\frac{p}{n}}\geq (|\bn n|\ |\bn n|)^{\frac{p}{n}}.
$$
As $n+p>0$,
\begin{equation*}\label{i-0}
\frac{as_p(K)}{as_p(B^n_2 )}\geq
\left(\frac{|K|}{|B^n_2|}\right)^{\frac{n-p}{n+p}}.
\end{equation*}
The equality case follows from  the equality case for the Blaschke Santal\'o
inequality.
\par
\vskip 3mm  (iii) Similarly, combining  (\ref{4.1.5c}),
$n|B^n_2|=(as_p(B^n_2))^{\frac{n+p}{p-1}}(as(B^n_2))^{\frac{n+1}{1-p}}$,
and the Inverse Santal\'o inequality, we get, for $p<-n$,
\begin{eqnarray*} \left(\frac{as_p(K)}{as_p(B^n_2
)}\right)^{\frac{n+p}{p}}\geq \left(\frac{|K^\circ|}{|B^n_2
|}\right)\left(\frac{|K|}{|B^n_2 |}\right)^{\frac{n}{p}} \geq c^n
\left(\frac{|K|}{|B^n_2 |}\right)^{\frac{n-p}{p}}.
\end{eqnarray*}
As $\frac{n+p}{p}>0$,
\begin{equation*} \label{i-1}
\frac{as_p(K)}{as_p(B^n_2 )}\geq
c^{\frac{np}{n+p}}\left(\frac{|K|}{|\bn
n|}\right)^{\frac{n-p}{n+p}}.
\end{equation*}

%---------------------------------------------------------------------------------
%-----------------------L_{-n} case-----------------------------------------------
%---------------------------------------------------------------------------------

\vskip 5mm The {\em  $L_{-n}$ affine surface area} was defined in
\cite{MW2} for convex bodies $K$ in $C^2_+$ and with centroid at
the origin by
$$
as_{-n}(K)=\max_{u\in S^{n-1}}f_K(u) ^{\frac{1}{2}}
h_K(u)^{\frac{n+1}{2}}.
$$

More generally, one could  define the   $L_{-n}$ affine surface
area for any convex body $K$ with centroid at the origin by
$
as_{-n}(K)=\sup _{x\in
\partial K} \frac{\langle x, N_K(x)\rangle^{\frac{n+1}{2}}}{\kappa
_K(x)^{\frac{1}{2}}}
$.
But as  in most cases then $as_{-n}(K) = \infty$, it suffices to consider   $K$ in $C^2_+$.
\vskip 4mm
A statement similar to Theorem \ref{p-aff-sio8} holds.
\vskip 3mm
\begin{proposition} \label{p-aff-sio9}

Let $K$ be a convex body in $C^2_+$ with centroid  at the
origin.  Let $p\neq -n$ and $s\neq -n$ be real
numbers.
\par
(i) If $\frac{n(s-p)}{(n+p)(n+s)}\geq0$, then

\begin{equation*} \label{T-4.3-1} as_p(K)\leq
\big(as_{-n}(K)\big)^{\frac{2n(s-p)}{(n+p)(n+s)}} \ as_s(K).
\end{equation*}
\par
(ii) If $\frac{n(s-p)}{(n+p)(n+s)}\leq 0$, then

\begin{equation*} \label{T-4.3-2} as_p(K)\geq
\big(as_{-n}(K)\big)^{\frac{2n(s-p)}{(n+p)(n+s)}} \ as_s(K).
\end{equation*}
\par
(iii) The  $L_{-n}$ affine isoperimetric inequality holds
\begin{equation*}\frac{as_{-n} (K)}{as_{-n}(\bn n)}\geq
\frac{|K|}{|\bn n|}.
\end{equation*}
\end{proposition}
\vskip 3mm
\noindent
{\bf Proof} \ \     (i) and (ii)
\begin{eqnarray*} as_p(K)&=&\int _{\partial K}
\frac{\kappa_K(x)^{\frac{p}{n+p}}}{\langle x, N_K(x) \rangle
^{\frac{n(p-1)}{n+p}}}\,d\mu _K (x)
\\ \\ &=&\int _{\partial K}
\left(\frac{\kappa_K(x)^{\frac{s}{n+s}}}{\langle x, N_K(x)
\rangle^{\frac{n(s-1)}{n+s}}}\right) \ \left( \frac{\langle x,
N_K(x)
\rangle^{\frac{n+1}{2}}}{\kappa_K(x)^{\frac{1}{2}}}\right)^{\frac{2n(s-p)}{(n+p)(n+s)}}\,d\mu_K(x)\end{eqnarray*}
which is
$$
\leq \big(as_{-n}(K)\big)^{\frac{2n(s-p)}{(n+p)(n+s)}} \ as_s(K),
\quad \mbox{if $\frac{n(s-p)}{(n+p)(n+s)} \geq 0$,}
$$
 and
$$
\geq \big(as_{-n}(K)\big)^{\frac{2n(s-p)}{(n+p)(n+s)}} \ as_s(K),
\quad \mbox{if $\frac{n(s-p)}{(n+p)(n+s)} \leq 0$.}
$$

\par (iii)
Note that
 $\frac{n(s-p)}{(n+p)(n+s)}>0$
implies  that $s>p>-n$ or $p<s<-n$ or $s<-n<p$. If $p=0$ and
$s\rightarrow \infty$, then
\begin{equation} \label{T-4.3-3}
as_{-n} (K) \geq \sqrt{\frac{|K|}{|K^\circ|}}.
\end{equation}

This gives the $L_{-n}$ affine isoperimetric inequality
\begin{equation*} \label{T-4.3-4}
\frac{as_{-n} (K)}{as_{-n}(\bn n)}=as_{-n}(K) \geq
\sqrt{\frac{|K|}{|K^\circ|}}\geq \sqrt{\frac{|K|^2}{|K|\
|K^\circ|}}\geq \sqrt{\frac{|K|^2}{|\bn n|^2}}=\frac{|K|}{|\bn
n|}.
\end{equation*}

\vskip 5mm

Analogous to corollary \ref{p-aff-iso3}, an immediate consequence
of (\ref{T-4.3-3}) is the following corollary. It can also be
proved directly using (\ref{curvatureformula}).
\vskip 3mm

\begin{corollary}
 Let
$K$ be a convex body in $C^2_+$ with centroid at the origin. Then
$$as_{-n}(K^\circ)as_{-n}(K) \geq {as_{-n}(\bn n)^2}.$$
\end{corollary}

\vskip 5mm
\noindent
{\bf Acknowledgment} The authors would like to thank the referee for the many  helpful suggestions.

\vskip 2mm \noindent Elisabeth Werner\\
{\small Department of Mathematics \ \ \ \ \ \ \ \ \ \ \ \ \ \ \ \ \ \ \ Universit\'{e} de Lille 1}\\
{\small Case Western Reserve University \ \ \ \ \ \ \ \ \ \ \ \ \ UFR de Math\'{e}matique }\\
{\small Cleveland, Ohio 44106, U. S. A. \ \ \ \ \ \ \ \ \ \ \ \ \ \ \ 59655 Villeneuve d'Ascq, France}\\
{\small \tt elisabeth.werner@case.edu}\\ \\
\and Deping Ye\\
{\small Department of Mathematics}\\
{\small Case Western Reserve University}\\
{\small Cleveland, Ohio 44106, U. S. A.}\\
{\small \tt dxy23@case.edu}
\end{document}